\documentclass[12pt,psamsfonts]{article}
\usepackage{amsmath}
\usepackage{amssymb}
\usepackage{amsthm}

\newtheorem{defn0}{Definition}[section]
\newtheorem{prop0}[defn0]{Proposition}
\newtheorem{thm0}[defn0]{Theorem}
\newtheorem{lemma0}[defn0]{Lemma}
\newtheorem{corollary0}[defn0]{Corollary}
\newtheorem{example0}[defn0]{Example}
\newtheorem{remark0}[defn0]{Remark}
\newtheorem{conjecture0}[defn0]{Conjecture}

\newenvironment{proposition}{\bigskip \begin{prop0}}{\end{prop0}}
\newenvironment{theorem}{\bigskip \begin{thm0}}{\end{thm0}}
\newenvironment{lemma}{\bigskip \begin{lemma0}}{\end{lemma0}}

\def\cocoa
{\mbox{\rm C\kern-.13em o\kern-.07
em C\kern-.13em o\kern-.15em A}}

\newcommand{\m}{\mathfrak{m}}
\newcommand{\n}{\mathfrak{n}}

\newif\ifprivate
\privatetrue

\title{  \bf \huge Isomorphism classes of certain Artinian Gorenstein Algebras.
}
\author{\large   Juan Elias
\thanks{Partially supported by  MTM2007-67493}
\and \large Giuseppe Valla \
\thanks{Partially supported by MIUR under the framework of PRIN 2007}}

\date{\today}

\begin{document}

\maketitle

\bigskip

\noindent
Juan Elias\\
Departament d'\`Algebra i Geometria\\
Universitat de Barcelona\\
Gran Via 585, 08007 Barcelona, Spain\\
e-mail: {\tt elias@ub.edu}

\bigskip
\noindent
Giuseppe Valla\\
Dipartimento di Matematica\\
Universit{\`a} di Genova\\
Via Dodecaneso 35, 16146 Genova, Italy\\
e-mail: {\tt valla@dima.unige.it}

\begin{abstract}  In this paper we attack the problem of the classification, up to analytic isomorphism, of Artinian Gorenstein local $k$-algebras with a given Hilbert Function.  We solve the problem in the case  the square of the maximal ideal is minimally generated by two elements and  the socle degree is high enough.
\end{abstract}

\bigskip
\section{Introduction.}

A longstanding problem in Commutative Algebra
is the classification of Artin algebras. We know that there exists a
finite number of isomorphism classes of Artin algebras of
multiplicity at most $6$,  however,  if the multiplicity is at least
$7$ there are examples with infinitely many isomorphism classes, see \cite{Iar87b}
and \cite{Poo07} and their reference lists for more results on the
classification problem. The aim of this paper is to classify the
family of almost stretched Gorenstein Artin algebras.

 In \cite{Sal79c} a local Artinian ring
$(A,\m)$ is said to be stretched if the square of the maximal  ideal $\m$ of $A$
is a principal ideal. In that paper J. Sally gave a nice  structure
theorem for stretched Gorenstein local rings.  Other interesting
properties of stretched $\m$-primary ideals can be found in
\cite{RosV01}. Sally's result has been considerably extended in
\cite{EV}, where   the notion of almost stretched local rings has
been introduced. A local Artinian ring $(A,\m)$ is said to be {\bf
almost stretched} if the minimal number of generators of $\m^2$ is
two.

We know from the classical Theorem of Macaulay, concerning the
possible Hilbert  functions of standard graded algebras, that the
Hilbert function of an almost stretched Gorenstein local ring $A$
has the following shape  $(1,h,2,...,2,1,...,1),$ which means
$$H_A(j)=
\begin{cases}
     1 & \ \ \text{$j=0$}, \\
     h & \ \ \text{$j=1$},\\
     2 & \ \ \text{$2\le j \le t$},\\
     1 & \ \ \text{$t+1\le j \le s,$}
\end{cases}$$  for integers $s$ and $t$ such that $s\ge t+1 \ge 3$ and $h\ge 2.$ If an algebra has  this Hilbert Function we say that it is {\bf of type $(s,t).$}

In \cite{EV} we gave a useful    structure theorem  for  almost
stretched Gorenstein local rings in the embedded  case, namely when
$A=R/I$ with $(R,\n)$ a regular local ring of dimension $h$ such
that $k:=R/\n$ is algebraically closed of characteristic 0.

The  result reads as follows. Let $I$ be an ideal
of $R$; then $A:=R/I$ is almost stretched and Gorenstein of type $(s,t)$
 if and only if there exists a minimal
basis $y_1,\dots,y_h$ of $\n$ and an element $b\in R$, such that
the ideal $I$ is generated by the following ${h \choose 2} +h-1$ elements:

$y_iy_j$, with $1\le i < j \le h,\ \ \ (i,j)\not=(1,2),$

$ y_j^2-y_1^s$ with $3\le j\le h,$

 $y_2^2-by_1y_2-y_1^{s-t-1},$

 $ y_1^ty_2.$

In this paper, we assume $R$ to be a power series ring of dimension $h$  over an algebraically closed field $k$ and  we fix  the  integers $s,t$ such that $s\ge t+1 \ge 3.$  We attack the problem of classifying, up to  analytic
isomorphism, the family of almost stretched  Gorenstein Algebras  $A=R/I$  of type $(s,t).$
   We solve the
problem in the case the socle degree is large enough with respect to $t,$ namely when $s\ge 2t,$ see Theorem \ref{finale}.

It turns out that for generic $s$ and $t$   we have exactly $t$ isomorphism classes of Gorenstein almost stretched algebras. But  if there  exists an integer $r$ with the properties  $0\le r\le t-2$ and $2(r+1)=s-t+1,$ then one dimensional families of non isomorphic models arise. Notice that, as a particular case,  we prove the existence of a codimension two complete intersection algebra of length 10 with infinitely many isomorphism classes.

Suitable examples at the end of the paper show that the case $s\le 2t-1$ is far away from a solution.

\section{The models.}

\label{preli} Through the paper  we are assuming that the basic
field $k$ is an algebraically  closed field of characteristic zero.

We will also freely use the following result which is a
straightforward application of Hensel  lemma.

\begin {proposition} \label{Hensel-like}
 Let $f=f(x_1,\dots,x_n)\in k[[x_1,\dots,x_n]]$ be  an
invertible  formal power series with $f(0,\dots,0)=a_0\not=0.$ If
there exists $ \alpha\in k$  such that $\alpha^j=a_0,$ then there
exists $g\in R$ such that $g^j=f$ and $g(0,\dots, 0)=\alpha$.
\end{proposition}

Let $R=k[[x_1,\dots, x_h]]$  be the formal power series ring and $ \n$ its maximal ideal. Given a  set of generators $\underline{y}=\{y_1,\dots,y_h\}$
of $\n,$ we let $\varphi_{\underline{y}}$ be the automorphism of $R$ which
is the result of  substituting $y_i$ for $x_i$ in a
power series $f(x_1,\dots,x_h)\in R.$  It is well known that given two ideals
$I$ and $J$ in $R$  there exists a $k$-algebras isomorphism
$$\alpha: R/I\to R/J$$ if and only if for some generators $\{y_1,\dots,y_h\}$ of
$\n$ we have $I=\varphi_{\underline{y}}(J).$

By abuse of notation, we will often say that $I$ is isomorphic to
$J$ and we will write  $I \sim J,$ with the meaning that $R/I$ is
isomorphic to $R/J.$

As we explained above, an ideal $I$ of $R$ is almost stretched and Gorenstein of type $(s,t)$ if and only if for suitable  formal power series $y_1,\dots,y_h,b\in R$  we have  $\n=(y_1,\dots,y_h)$ and $$I=(\underset{
\underset{(i,j)\not=(1,2)}
{1\le i < j \le h}}
{y_iy_j}, \underset{3\le j \le h}{y_j^2-y_1^s}, y_2^2-by_1y_2-y_1^{s-t+1},y_1^ty_2).$$ If we let $a:=a(x_1,\dots,x_h)$ be the formal power series such that $a(y_1,\dots,y_h)=b,$ we get  $I=\varphi_{\underline y}(I_a)$ where $$I_a:=(\underset{
\underset{(i,j)\not=(1,2)}
{1\le i < j \le h}}
{x_ix_j}, \underset{3\le j \le h}{x_j^2-x_1^s}, x_2^2-ax_1x_2-x_1^{s-t+1},x_1^tx_2).$$

Hence we can  rephrase  the main result in \cite{EV} as follows.
An ideal $I$ of $R$ is almost stretched and Gorenstein of type $(s,t)$ if and only if it is isomorphic to $I_a$ for some $a\in R.$

This implies that, in order to get a  classification of the  almost stretched  Gorenstein Algebras  $A=R/I$  of type $(s,t),$ we need to classify the algebras $A=R/I_a$ when $a$ is running in $R.$

We will see that the order of the power series $a(x_1,0,\dots,0)$ plays a
central role in the  classification problem. Hence, first we study
the case $a(x_1,0,\dots,0)=0.$

 Given the integer $p\ge 0$ and the power series $z\in R,$ we introduce the ideal $I_{p,z}$ which is  generated as follows:
$$I_{p,z}:=(\underset{
\underset{(i,j)\not=(1,2)}
{1\le i < j \le h}}
{x_ix_j}, \underset{3\le j \le h}{x_j^2-x_1^s}, x_2^2-x_1^{p+1}x_2-zx_1^{s-t+1},x_1^tx_2)$$
These ideals will be crucial in the rest of the paper.

\begin{proposition}\label{a(x,0)=0} If $a(x_1,0,\dots,0)=0,$ then $$I_a\sim I_{t-1,1}.$$
\end{proposition}
\begin{proof} First we remark that $$I_{t-1,1}=
(\underset{
\underset{(i,j)\not=(1,2)}
{1\le i < j \le h}}
{x_ix_j}, \underset{3\le j \le h}{x_j^2-x_1^s}, x_2^2-x_1^tx_2-x_1^{s-t+1},x_1^tx_2)=(\underset{
\underset{(i,j)\not=(1,2)}
{1\le i < j \le h}}
{x_ix_j}, \underset{3\le j \le h}{x_j^2-x_1^s}, x_2^2-x_1^{s-t+1},x_1^tx_2).$$
If $a(x_1,0,\dots,0)=0,$ we have  $a=\sum_{i=2}^hb_ix_i$ with $b_i\in R,$ so that
$$x_2^2-ax_1x_2-x_1^{s-t+1}=x_2^2-x_1x_2(\sum_{i=2}^hb_ix_i)-x_1^{s-t+1}=x_2^2(1-b_2x_1)-x_1x_2(\sum_{i=3}^hb_ix_i)-x_1^{s-t+1}.$$
By Proposition \ref{Hensel-like} we can find a power series $v\in R$ such that $v^2=1-b_2x_1.$
Then $v\notin \n$ and we have
$$I_a:=(\underset{
\underset{(i,j)\not=(1,2)}
{1\le i < j \le h}}
{x_ix_j}, \underset{3\le j \le h}{x_j^2-x_1^s}, (x_2v)^2-x_1^{s-t+1},x_1^t(x_2v))=$$
$$=(\underset{3\le j\le h}{x_2vx_j},\underset{
\underset{(i,j)\not=(1,2), i\not=2}
{1\le i < j \le h}}
{x_ix_j}, \underset{3\le j \le h}{x_j^2-x_1^s}, (x_2v)^2-x_1^{s-t+1},x_1^t(x_2v)).$$ If we consider the change of variables $$y_j:=\begin{cases}
  x_j    & \text{if $j\not=2$} \\
    vx_2  & \text{if $j=2$}
\end{cases}$$  then we have $\n=(y_1,\dots,y_h)$ and $$I_a=(\underset{
\underset{(i,j)\not=(1,2)}
{1\le i < j \le h}}
{y_iy_j}, \underset{3\le j \le h}{y_j^2-y_1^s}, y_2^2-y_1^{s-t+1},y_1^ty_2)=\varphi_{\underline{y}}(I_{t-1,1}).$$
The conclusion follows.
 \end{proof}

    Next we study the case  $a(x_1,0,\dots,0)\not=0.$
\begin{proposition}  \label{Irw}
 If $a(x_1,0,\dots,0)\not=0,$ then  $I_a\sim I_{r,w},$  where $w$ is a suitable  invertible power series in $R$ and $r$ is the order of  $a(x_1,0,\dots,0)$ in $k[[x_1]].$ Further, if $s\ge 2t-1,$ we may assume $w\in k[[x_1]]\setminus (x_1).$
 \end{proposition}
 \begin{proof} Since $r$ is the order of $a(x_1,0,\dots,0)$ in $k[[x_1]],$ we can write $a=x_1^r\eta+\sum_{j=2}^hb_jx_j$ with $b_j\in R,$ $\eta \in k[[x_1]]\setminus (x_1).$ We get
 $$x_2^2-ax_1x_2-x_1^{s-t+1}=x_2^2-x_1x_2( x_1^r\eta+\sum_{j=2}^hb_jx_j)-x_1^{s-t+1}=$$
 $$=x_2^2-x_1^{r+1}x_2\eta-x_1x_2(\sum_{j=2}^hb_jx_j)-x_1^{s-t+1}=x_2^2(1-b_2x_1)-x_1^{r+1}x_2\eta-\sum_{j=3}^hb_jx_jx_1x_2-x_1^{s-t+1}.$$
 If we let $u:=1-x_1b_2,$ then $u\notin \n$ and,  by Proposition \ref{Hensel-like}, we can find a power series $z\notin \n$ such that $z^{r+1}=\eta/u.$ Then we get
 $$I_a=(\underset{
\underset{(i,j)\not=(1,2)}
{1\le i < j \le h}}
{x_ix_j}, \underset{3\le j \le h}{x_j^2-x_1^s}, ux_2^2-x_1^{r+1}x_2\eta-x_1^{s-t+1},x_1^tx_2)=$$
$$=(\underset{
\underset{(i,j)\not=(1,2)}
{1\le i < j \le h}}
{x_ix_j}, \underset{3\le j \le h}{x_j^2-x_1^s},x_2^2-x_1^{r+1}x_2z^{r+1}-\frac{(x_1)^{s-t+1}}{u},x_1^tx_2)=$$
$$=(\underset{
\underset{(i,j)\not=(1,2)}
{1\le i < j \le h}}
{x_ix_j}, \underset{3\le j \le h}{x_j^2-x_1^s},x_2^2-(x_1z)^{r+1}x_2-\frac{(x_1z)^{s-t+1}}{uz^{s-t+1}},x_1^tx_2).$$
By Proposition  \ref{Hensel-like} we can find a power series $\tau\in R\setminus \n$ such that $\tau^2=z^s$; now we consider  the following change of variables $$y_j:=\begin{cases}
   zx_1   & \text{if $j=1$}, \\
  x_2    & \text{if $j=2$}\\
  \tau x_j & \text{if $j=3,\dots,h$}.
\end{cases}$$ It is clear that $\n=(y_1,\dots,y_h);$  further $z^s(x_j^2-x_1^s)=\tau^2x_j^2-z^sx_1^s=y_j^2-y_1^s$
and we can find an invertible  power series $w:=w(x_1,\dots,x_h)$ such that $w(y_1,\dots,y_h)=\frac {1}{uz^{s-t+1}}.$
Hence we get
$$I_a=(
\underset{j=3,\dots,h}{(zx_1)(\tau x_j)},
\underset{j=3,\dots,h}{x_2(\tau x_j)},
\underset{3\le i < j \le h}{(\tau x_i)(\tau x_j)},
\underset{3\le j \le h}{z^s(x_j^2-x_1^s)},
x_2^2-(x_1z)^{r+1}x_2-\frac{(x_1z)^{s-t+1}}{uz^{s-t+1}},(zx_1)^tx_2
)=$$ $$=(\underset{\underset{(i,j)\not=(1,2)}{1\le i < j \le h}}{y_iy_j},
\underset{3\le j \le h}{y_j^2-y_1^s},
y_2^2-y_1^{r+1}y_2-w(y_1,\dots,y_h)y_1^{s-t+1},y_1^ty_2)=\varphi_{\underline{y}}(I_{r,w}).$$
This proves $I_a\sim I_{r,w}$ and the first assertion.

  As for the second one, let $s\ge 2t-1.$ It is clear that $wx_1^{s-t+1}=[w(x_1,0,\dots,0)+\sum_{j\ge 2}c_jx_j]x_1^{s-t+1}=w(x_1,0,\dots,0)x_1^{s-t+1}+\rho$ where $\rho\in (x_1x_3,\dots,x_1x_h,x_1^{s-t+1}x_2).$ Since $s-t+1\ge t,$ we have $x_1^{s-t+1}x_2\in (x_1^tx_2).$ By replacing $w$ with $w(x_1,0,\dots,0)$ we get  the conclusion.
 \end{proof}

 \vskip 2mm From the above result we need to study the isomophism classes of the ideals $I_{r,w}$ where {\bf $r$ is a non negative integer and $w$ an invertible power series}.
 \begin{proposition} \label{r ge t-1, s-t}
 If $r\ge t-1$ then $$
 I_{r,w}\sim  I_{t-1,1}.$$
\end{proposition}
\begin{proof}  If  $r\ge t-1$, then $r+1\ge t$ and we have
$$I_{r,w}=(\underset{\underset{(i,j)\not=(1,2)}{1\le i < j \le h}}{x_ix_j},
\underset{3\le j \le h}{x_j^2-x_1^s},
x_2^2-wx_1^{s-t+1},x_1^tx_2).$$ Let $v^2=\frac{1}{w},$ so that $v\notin \n$ and
$$I_{r,w}=(\underset{\underset{(i,j)\not=(1,2)}{1\le i < j \le h}}{x_ix_j},
\underset{3\le j \le h}{x_j^2-x_1^s},
w(\frac{x_2^2}{w}-x_1^{s-t+1}),x_1^tx_2)=$$
$$=(\underset{\underset{(i,j)\not=(1,2)}{1\le i < j \le h}}{x_ix_j},
\underset{3\le j \le h}{x_j^2-x_1^s},
(x_2v)^2-x_1^{s-t+1},x_1^tx_2)=$$
$$=(\underset{3\le j \le h}{x_1x_j},\underset{3\le j \le h}{(x_2v)x_j},\underset{3\le i<j\le h} {x_ix_j},
\underset{3\le j \le h}{x_j^2-x_1^s}, (x_2v)^2-x_1^{s-t+1},x_1^t(x_2v)).$$
We consider  the following change of variables $$y_j:=\begin{cases}
   x_j   & \text{if $j\not= 2$}, \\
  vx_j    & \text{if $j=2$}.
\end{cases}$$ It is clear that $\n=(y_1,\dots,y_h)$ and further
$$ I_{r,w}=(\underset{\underset{(i,j)\not=(1,2)}{1\le i < j \le h}}{y_iy_j},
\underset{3\le j \le h}{y_j^2-y_1^s},
y_2^2-y_1^{s-t+1},y_1^ty_2)=\varphi_{\underline{y}}(I_{t-1,1}).$$
This proves the result. \end{proof}

We can also understand the isomorphism class of the ideal $I_{r,w}$ in the case $r$ is general enough. More precisely we have the following result.

\begin{proposition} \label{r good} Let $r\ge 0$ be an integer such that $2(r+1)\not= s-t+1.$ Then $I_{r,w}\sim I_{r,1}$.
\end{proposition}
\begin{proof}  Let $n:=2(r+1)-( s-t+1)$  and choose a power series $e\in R$ such that $$\begin{cases}
  e^n  = \frac{1}{w} & \text{if $n>0$} \\
   e^{-n}=w   & \text{if  $n<0$}.
\end{cases}$$ In both cases we have $e^nw=1.$ Further let us choose an element $\tau \in R$ such that $\tau ^2=e^s.$ It is clear that both $e$ and $\tau $ are not in $ \n.$  We have
$$I_{r,w}=(\underset{\underset{(i,j)\not=(1,2)}{1\le i < j \le h}}{x_ix_j},
\underset{3\le j \le h}{x_j^2-x_1^s},
x_2^2-x_1^{r+1}x_2-wx_1^{s-t+1},x_1^tx_2)=$$
$$=(\underset{3\le j \le h}{ex_1(\tau x_j)},
\underset{3\le j \le h}{(e^{r+1}x_2)(\tau x_j)},
\underset{3\le i<j \le h}{\tau^2 x_ix_j},
\underset{3\le i<j \le h}{\tau^2 (x_j^2-x_1^s)},
e^{2r+1}(x_2^2-x_1^{r+1}x_2-wx_1^{s-t+1}), e^{t+r+1}x_1^tx_2),
$$ Now we remark that $$\tau^2 (x_j^2-x_1^s)=(\tau x_j)^2-(ex_1)^s,\ \ \ \ \ e^{t+r+1}x_1^tx_2=(ex_1)^t(e^{r+1}x_2)$$
$$(e^{r+1}x_2)^2-(ex_1)^{r+1}(e^{r+1}x_2)-(ex_1)^{s-t+1}=e^{2(r+1)}x_2^2-e ^{2(r+1)}x_1^{r+1}x_2-e^{s-t+1}x_1^{s-t+1}=
$$
$$=e^{2(r+1)}(x_2^2-x_1^{r+1}x_2-wx_1^{s-t+1})+we^{2(r+1)}x_1^{s-t+1}-e^{s-t+1}x_1^{s-t+1}=$$
$$=e^{2(r+1)}(x_2^2-x_1^{r+1}x_2-wx_1^{s-t+1})+x_1^{s-t+1}(we^{2(r+1)}-e^{s-t+1})=$$
$$=e^{2(r+1)}(x_2^2-x_1^{r+1}x_2-wx_1^{s-t+1})+e^{s-t+1}x_1^{s-t+1}(we^n-1)=$$
$$=e^{2(r+1)}(x_2^2-x_1^{r+1}x_2-wx_1^{s-t+1}).$$
We let $$y_j:=\begin{cases}
   ex_1   & \text{if $j=1$}, \\
  e^{r+1}x_2    & \text{if $j=2$}\\
  \tau x_j & \text{if $j=3,\dots,h$}.
\end{cases}$$ It is then clear that $\n=(y_1,\dots,y_h)$ and

$$I_{r,w}=(\underset{\underset{(i,j)\not=(1,2)}{1\le i < j \le h}}{y_iy_j},
\underset{3\le j \le h}{y_j^2-y_1^s},
y_2^2-y_1^{r+1}y_2-y_1^{s-t+1},y_1^ty_2)=\varphi_{\underline{y}}(I_{r,1}).$$ The result is proved.

\end{proof}

From the above result it becomes relevant the following notation. We say that the couple $(s,t)$ is {\bf regular} if  it does not exist an integer $r$ with the properties  $0\le r\le t-2$ and $2(r+1)=s-t+1$. It is easy to see that the couple $(s,t)$ is not regular if and only if $s-t$ is odd and $s\le 3t-3.$

For a regular couple $(s,t)$  we have at most $t$ isomorphism classes with models the ideals $I_{0,1},I_{1,1},\dots,I_{t-1,1}.$

Unfortunately it is not true, even in the case of two variables, that for a regular couple $(s,t)$ the above models are pairwise non isomorphic. At the end of the paper we will show that if  $h=2$ and we consider the couple $s=5,$ $t=3,$  which corresponds to the  Hilbert Function  $1,2,2,1,1,1,$ then the couple $(5,3)$ is regular but  $$I_{1,1}=(y^2-x^2y-x^3,x^3y)\simeq I_{2,1}=(y^2-x^3,x^3y).$$

In the above example we have $s< 2t.$ Namely we will prove that, if $s\ge 2t,$ then the ideals  $I_{0,1},I_{1,1},\dots,I_{t-1,1}$ are pairwise non isomorphic.

%We start by comparing $I_{0,1}$ with $I_{j,1},$ $j\ge 1.$
%\begin{proposition} \label{r=0} If $s\ge 2t$ then $I_{0,1}$ is not isomorphic to $I_{j,1},$ for every $j\ge 1.$
%\begin{proof} We first remark that $$(x_1-x_2)(x_1^{s-t}+x_2)=-(x_2^2-x_1x_2-x_1^{s-t+1})-x_1^{s-2t}(x_1^tx_2)\in I_{0,1}.$$ If we assume by contradiction that $I_{j,1}=\phi_{\underline y}(I_{j,1}),$ with $\n=(y_1,\dots,y_h),$ then
%$$\phi_{\underline y}((x_1-x_2)(x_1^{s-t}+x_2))=(y_1-y_2)(y_1^{s-t}+y_2)\in I_{j,1}.$$ We let $$z_j:=\begin{cases}
 %  y_1-y_2   & \text{if $j=1$}, \\
  %y_1^{s-t}+y_2    & \text{if $j=2$}\\
 % y_j & \text{if $j=3,\dots,h$}.
%\end{cases}.$$ Since $$det\begin{pmatrix}
 % 1    & y_1^{s-t-1} &0 & \dots & 0    \\
  % -1   &  1 & 0  & \dots & 0 \\
 %  0    &  0 & 1 &  \dots & 0 \\
   %\vdots & \vdots & \vdots &\ddots & \vdots \\
   %0 & 0 & 0 & \dots &1
  %\end{pmatrix}= 1+x_1^{s-t-1}\notin \n,$$ we have

%\end{proof}\end{proposition}

  \vskip 3mm We need now to study the ideals $I_{r,w},$ where   $2(r+1)=s-t+1$ and $w$ is an invertible power series. We recall that, by  Proposition \ref{Irw}, when $s\ge 2t-1,$ we can further assume that $w\in k[[x_1]]\setminus (x_1).$

\begin{proposition} Let  $r\ge 0$ be an integer such that   $2(r+1)=s-t+1.$ Further let $w=\sum_{i\ge 0} w_ix_1^i$ be an invertible power series in $k[[x_1]]$ such that $w\not= w_0.$  If $d$ is the order of $w-w_0,$ then $I_{r,w}\sim I_{r,w_0+x_1^d}.$
\end{proposition}
\begin{proof} We have $w=w_0+w_dx_1^d+\dots$ with $w_d\in k^*,$ so that  we can find  a power series $\alpha \in k[[x_1]]$ such that $\alpha^d=\sum_{n\ge d}w_nx_1^{n-d}.$ Hence  $\alpha$ is invertible and 
$\alpha^d x_1^d=w-w_0.$ As a consequence we can find an invertible power series $\beta \in k[[x_1]]$ such that $\beta^2=\alpha ^s.$  Let us consider the following change of variables $$y_j:=\begin{cases}
   \alpha x_1   & \text{if $j=1$}, \\
  \alpha^{r+1}x_2    & \text{if $j=2$}\\
  \beta  x_j & \text{if $j=3,\dots,h$}.
\end{cases}$$

For every $j\ge 3$ we have $$\alpha ^s(x_j^2-x_1^s)=(\beta x_j)^2-(\alpha x_1)^s=y_j^2-y_1^s;$$ further
$$\alpha^{t+r+1}x_1^tx_2=(\alpha x_1)^t(\alpha^{r+1}x_2)=y_1^ty_2$$ and
$$y_2^2-y_1^{r+1}y_2-(w_0+y_1^d)y_1^{s-t+1}=
x_2^2\alpha^{2(r+1)}-x_1^{r+1}\alpha^{r+1}x_2\alpha^{r+1}-(w_0+y_1^d)x_1^{2(r+1)}\alpha^{2(r+1)}=
$$
$$=\alpha^{2(r+1)}(x_2^2-x_2x_1^{r+1}-wx_1^{2(r+1)})
$$ where the last equality follows because 
$$
y_1^d=\alpha^d x_1^d=w-w_0.
$$
As a consequence we get that the ideal $I_{r,w}$ coincides with the ideal
$$(
\underset{3\le j \le h}{\alpha x_1\beta x_j},
\underset{3\le j \le h}{\alpha^{r+1}x_2\beta x_j},
\underset{3\le i < p\le h}{\beta x_i\beta x_p},
\underset{3\le j \le h}{\alpha^s(x_j^2-x_1^s)}, \alpha^{2(r+1)}(x_2^2-x_2x_1^{r+1}-wx_1^{2(r+1)}),\alpha^{t+r+1}x_1^tx_2)=$$
$$=(\underset{3\le j \le h}{y_1y_j},
\underset{3\le j \le h}{y_2y_j},
\underset{3\le i < p\le h}{y_iy_p},
\underset{3\le j \le h}{y_j^2-y_1^s}, y_2^2-y_2y_1^{r+1}-(w_0+y_1^d)y_1^{2(r+1)},y_1^ty_2)=\varphi_{\underline{y}}(I_{r,w_0+x_1^d}).$$ This gives the conclusion.
\end{proof} We can now improve the last result in the case $d\ge t-r-1.$

\begin{proposition} Let $r\ge 0$ and $d$ an integer such that $d\ge t-r-1.$ Then for every $c \in k^*$ we have $I_{r,c +x_1^d}\sim I_{r,c}.$
\end{proposition}
\begin{proof} For simplicity we let $\eta:=c+x_1^d.$ Let $\alpha$ be a power series such that $\alpha^2=\frac{c}{\eta}.$  We must have $\alpha_0^2=1$ so that we may choose $\alpha$ with the properties $\eta\alpha^2=c$ and $\alpha_0=1.$  We have $$(\alpha x_2)^2-(\alpha x_2)x_1^{r+1}-cx_1^{s-t+1}=$$
$$=\alpha^2(x_2^2-x_1^{r+1}x_2-\eta x_1^{s-t+1})+\alpha^2x_2x_1^{r+1}+\alpha^2 \eta x_1^{s-t+1}-\alpha x_2x_1^{r+1}-cx_1^{s-t+1}=$$  $$=\alpha^2(x_2^2-x_1^{r+1}x_2-\eta x_1^{s-t+1})+\alpha x_2x_1^{r+1}(\alpha-1).$$  Now we have $$(\alpha-1)(\alpha+1)=\alpha^2-1=\frac{c}{\eta}-1=-\frac{x_1^d}{\eta}
\in (x_1^{t-r-1})$$ because $d\ge  t-r-1.$ Further, since $\alpha_0=1,$ the power series $\alpha +1$ is invertible. This implies $\alpha -1\in (x_1^{t-r-1}).$ It follows that for  suitable $\beta \in R$ we have $$(\alpha x_2)^2-(\alpha x_2)x_1^{r+1}-cx_1^{s-t+1}=\alpha^2(x_2^2-x_1^{r+1}x_2-\eta x_1^{s-t+1})+\beta x_1^tx_2.$$ Thus we get
$$I_{r,\eta}=(
\underset{3\le j \le h}{x_1x_j},
\underset{3\le j \le h}{\alpha x_2x_j},
\underset{3\le i < j \le h}{x_ix_j},
\underset{3\le j \le h}{x_j^2-x_1^s},
\alpha^2(x_2^2-x_1^{r+1}x_2-\eta x_1^{s-t+1}),x_1^t(\alpha x_2)
)=$$
$$=(
\underset{3\le j \le h}{x_1x_j},
\underset{3\le j \le h}{\alpha x_2x_j},
\underset{3\le i < j \le h}{x_ix_j},
\underset{3\le j \le h}{x_j^2-x_1^s},
(\alpha x_2)^2-(\alpha x_2)x_1^{r+1}-cx_1^{s-t+1},x_1^t(\alpha x_2)
).$$ We change the variables as follows  $$y_j:=\begin{cases}
    x_j  & \text{if $j\not= 2$}, \\
  \alpha x_2    & \text{if $j=2$}
  \end{cases}$$ and we get
  $$I_{r,\eta}=(\underset{\underset{(i,j)\not=(1,2)}{1\le i < j \le h}}{y_iy_j},
\underset{3\le j \le h}{y_j^2-y_1^s},
y_2^2-y_1^{r+1}y_2-c y_1^{s-t+1},y_1^ty_2)
=\varphi_{\underline{y}}(I_{r,c}).$$ The conclusion follows.
\end{proof}

Collecting the results of this section, we have  the following Theorem.
\begin{theorem}\label{modelli}  Let $I$ be an ideal in $R=k[[x_1,\dots,x_h]]$ such that $R/I$ is almost stretched and Gorenstein of type $(s,t)$ with $s\ge 2t-1.$

If $(s,t)$ is regular, then $I$ is isomorphic to one of the following ideals:
$$I_{0,1},\  I_{1,1},\dots , I_{t-1,1}$$

If  $(s,t)$ is not regular and $r$ is the integer such that $2(r+1)=s-t+1,$ then $I$ is isomorphic to one of the following ideals:

  $$I_{0,1},\dots, I_{r-1,1},\{I_{r,c}\}_{c\in k^*},\{I_{r,c+x_1}\}_{c\in k^*},\dots,\{I_{r,c+x_1^{t-r-2}}\}_{c\in k^*},I_{r+1,1},\dots , I_{t-1,1},$$ 
  
  \noindent with the meaning  that, if $r=t-2,$ then  the list of the possible models is
$$I_{0,1},\dots, I_{r-1,1},\{I_{r,c}\}_{c\in k^*},I_{r+1,1},\dots , I_{t-1,1}.$$
\end{theorem}

\section{Non isomorphic models.}

We are going now to prove that if $s\ge 2t$ and however we choose $w_0,\dots,w_{t-1}$  in $R\setminus \n,$ two ideals in the following list are never isomorphic: $$I_{0,w_0},I_{1,w_1},\dots,I_{t-1,w_{t-1}}.$$

We will need frequently in this section the following result proved in  \cite{EV}, Lemma 4.5.

Let $v_1,\dots,v_h$  be a minimal system  of generators  of the maximal ideal $\n$ of $R$ and let  $I$ be the ideal
$$I:=(\underset{\underset{(i,j)\not=(1,2)}{1\le i < j \le h}}{v_iv_j},
\underset{3\le j \le h}{v_j^2-v_1^s},
v_2^2-av_1v_2-uv_1^{s-t+1},v_1^tv_2)$$ with   $a\in R$ and $u\in R\setminus \n.$
The elements $\overline{v_1}^j,\overline{v_1}^{j-1}\overline{v_2}\in R/I$  form a minimal basis of the ideal $(\n/I)^j$ in the case $2\le j\le t$, while in the case  $t+1\le j\le s$ the ideal $(\n/I)^j$ is a principal ideal generated by $\overline{v_1}^j.$ Further $\n^{s+1}\subseteq I.$

Of course this means that \begin{equation}\label{base} \begin{cases}
  av_1^j+b{v_1}^{j-1}v_2\in I+\n^{j+1} \Longrightarrow  a,b\in \n    & \text{if $2\le j \le t$ }, \\
  av_1^j\in I+\n^{j+1}  \Longrightarrow  a\in \n   & \text{if $t+1\le j\le s$}.
\end{cases}\end{equation}

Because of its relevance, it is perhaps  useful to recall also the following  notation already introduced in Section 2. If  $p\ge 0$ is an integer and $z$ an invertible  power series, we let $$I_{p,z}=(\underset{\underset{(i,j)\not=(1,2)}{1\le i < j \le h}}{x_ix_j},
\underset{3\le j \le h}{x_j^2-x_1^s},
x_2^2-x_1^{p+1}x_2-zx_1^{s-t+1},x_1^tx_2).$$

\begin{lemma}\label{l2} Let $I:=I_{p,z}$ and  $l:=\sum_{i=1}^ha_ix_i,$ with $a_i\in R.$ If $s\ge t+2$ and $l^2\in I+\n^3,$ then $a_1\in \n.$ If $p=0,$ then $a_1,a_2\in \n.$
\end{lemma}
\begin{proof}
 We denote by $\cong$ congruence modulo the ideal $I+\n^3.$ We have $$0\cong l^2= \sum_{j=1}^ha_j^2x_j^2+2 \sum_{1\le i<j\le h}a_ia_jx_ix_j\cong \sum_{j=1}^ha_j^2x_j^2+2a_1a_2x_1x_2.$$ Since $s\ge t+2,$ we get $s-t+1\ge 3$ so that $x_2^2\cong x_1^{p+1}x_2.$ Further, for $j\ge 3$ we have $x_j^2\cong 0,$ because $x_j^2\cong x_1^s$ and   $s\ge 3.$

Hence $$a_1^2x_1^2+a_2^2x_1^{p+1}x_2+2a_1a_2x_1x_2=a_1^2x_1^2+x_1x_2(a_2^2x_1^p+2a_1a_2)\cong 0.$$
Now, if $p>0,$ then  $a_1^2x_1^2+2a_1a_2x_1x_2\in I+\n^3;$ this implies $a_1^2\in \n$ and finally $a_1\in \n.$

 If instead  $p=0,$ then $$a_1^2x_1^2+x_1x_2(a_2^2+2a_1a_2)\in I+\n^3$$ which implies $a_1^2,a_2^2+2a_1a_2\in \n$ and thus $a_1,a_2 \in \n.$

\end{proof}

As a consequence of this lemma, we can prove that  in the case $s\ge t+2$ and however we choose the units $z$ and $w$  in $R,$ the ideal $I_{0,w}$ is not isomorphic to $I_{p,z}$ when $p>0.$

\begin{proposition}\label{I0} Let $s\ge t+2$ and $p\ge 1.$ However we choose  the units $z$ and $w$ in $R$, we have $I_{0,w}\not\sim I_{p,z}.$
\end{proposition}
\begin{proof} By contradiction let us assume that $I_{0,w}=\varphi_{\underline{y}}(I_{p,z}),$ where $(y_1,\dots,y_h)=\n.$ This means that $$I:=I_{0,w}=(\underset{\underset{(i,j)\not=(1,2)}{1\le i < j \le h}}{y_iy_j},
\underset{3\le j \le h}{y_j^2-y_1^s},
y_2^2-y_1^{p+1}y_2-z( \underline{y})y_1^{s-t+1},y_1^ty_2)
.$$ We have $p\ge 1$ and $s\ge t+2$ so that $p+2\ge 3$ and $s-t+1\ge 3.$ This implies $y_2^2\in I+\n^3.$
Further, for every $j\ge 3,$ $y_j^2\in I+\n^s\subseteq I+\n^3,$ because  $s\ge 3.$ By the above Lemma, we get for every $j\ge 2,$ $$y_j=\sum_{i=3}^hc_{ji}x_i+d_j$$ where $d_j\in \n^2.$ Then we get $$I+\n^3=(\underset{\underset{(i,j)\not=(1,2)}{1\le i < j \le h}}{y_iy_j},\underset{2\le j \le h}{y_j^2})+\n^3\subseteq(y_1y_3,\dots,y_1y_h,\underset{3\le i\le  j \le h}{x_ix_j})+\n^3.$$ Let us  consider the $R/\n$-vector spaces $$V:=(I+\n^3)/\n^3\ \ \ \ \ \ \  W:=[(y_1y_3,\dots,y_1y_h,\underset{3\le i\le  j \le h}{x_ix_j})+\n^3]/ \n^3.$$
We have $$\dim_{R/\n}V\le \dim_{R/\n}W\le h-2+\binom{h-2+2-1}{2}=\binom{h}{2}-1.$$ 
On the other hand $\n^3\subseteq I+\n^3 \subseteq\n^2\subseteq \n \subseteq R$ so that $\dim (I+\n^3)/\n^3=\binom{h+1}{2}-2.$ Since $h\ge 2$ we get a contradiction.
\end{proof}

We are going now to prove that if $s\ge 2t$ the ideals $I_{1,w_1},I_{2,w_2},\dots,I_{t-1,w_{t-1}}$ are pairwise non isomorphic. In order to do that, we need some preparatory result.
\begin{lemma}\label{monomi}  Every monomial of degree $t+1$ in the variables $x_1,\dots,x_h$ is in $I_{p,z}+\n^s$, except possibly for $x_1^{t+1}.$
\end{lemma}

\begin{proof} Let $I:=I_{p,z}.$ Since $x_ix_j\in I$ if $ i<j$ and $(i,j)\not= (1,2),$ and $x_j^2\in I+\n^s$ if $3\le j \le h,$ we need only to consider monomials in $x_1,x_2.$ Now we have $x_1^tx_2\in I$ and then we can use descending induction to prove the result. So let $0\le j\le t-1$ and $x_1^{j+1}x_2^{t-j}\in I+\n^s;$ we need to prove that $x_1^jx_2^{t+1-j}\in I+\n^s.$
We denote by $\cong$ congruence modulo the ideal $I+\n^s.$ We have

\begin{equation*}\begin{split}x_1^jx_2^{t+1-j}&=x_1^jx_2^{t-1-j}x_2^2\\ &\cong x_1^jx_2^{t-1-j}(x_1^{p+1}x_2+zx_1^{s-t+1})\\ & =x_1^{j+p+1}x_2^{t-j}+zx_1^{s-t+j+1}x_2^{t-j-1} \\ &\cong 0\end{split}
\end{equation*} because $x_1^{s-t+j+1}x_2^{t-j-1}\in \n^s$ and, by induction, $x_1^{j+1}x_2^{t-j}\in I+\n^s.$ The conclusion follows.
\end{proof}

\begin{lemma}\label{y2} Let $s\ge t+2,$  $r\ge 1$ and $I_{p,z}=\varphi_{\underline{y}}(I_{r,w})$ where $(y_1,\dots,y_h)=\n.$ Then we can write $y_2=cx_1^{s-t}+d,$ with $c\in R$ and $d\in (x_2,\dots,x_h).$
\end{lemma}
\begin{proof} We let $I:=I_{p,z};$ hence $$I=(\underset{\underset{(i,j)\not=(1,2)}{1\le i < j \le h}}{y_iy_j},
\underset{3\le j \le h}{y_j^2-y_1^s},
y_2^2-y_1^{r+1}y_2-w( \underline{y})y_1^{s-t+1},y_1^ty_2).$$ For every $j\ge 3$ we have $y_j^2\in I+\n^s\subseteq I+\n^3,$ because $s\ge 3.$ Further $y_2^2\in I+\n^{r+2}+\n^{s-t+1}\subseteq I+\n^3,$ because $r\ge 1$ and $s\ge t+2.$ By Lemma \ref{l2} we have   $y_j=\sum_{i=1}^ha_{ji}x_i$ with $a_{j1}\in \n$ for every $j\ge 2.$

Since   $y_1,\dots,y_h$ is a minimal system of generators for $\n,$ we  must  have $y_1=ex_1+f$ with $e\notin \n$ and $f\in (x_2,\dots,x_h).$  We can also write  $y_2=ax_1^2+b$ with $a\in R$  and $b\in (x_2,\dots,x_h).$   Now, if $s=t+2$ we are done. Let $s\ge t+3$ and by induction let $y_2=ax_1^j+b$ with $b\in (x_2,\dots,x_h)$ and $ 2\le j \le s-t-1.$ We  claim that $a\in \n$ and remark that this would imply  the lemma.

We have $y_1^ty_2=(ex_1+f)^t(ax_1^j+b)\in I.$ By the above lemma we get $e^tax_1^{t+j}\in I+\n^s$ and since $e$ is a unit,  $ax_1^{t+j}\in I+\n^s.$  Now, if also $a$ is  a unit, we would get   $x_1^{t+j}\in I+\n^s$ so that
$x_1^{t+j+1}\in I.$ Since $s\ge t+j+1,$ this implies $x_1^s\in I$, a contradiction.
\end{proof}

We can prove now the main result of this section.
\begin{theorem}\label{noniso} Let $s\ge 2t$ and $1\le r \le t-2.$ If $p>r$ and however we choose $z,w\notin \n,$ the ideals $I_{p,z}$ and $I_{r,w}$ are never isomorphic.
\end{theorem}\begin{proof} Let us assume by contradiction that $I:=I_{p,z}=\varphi_{\underline{y}}(I_{r,w})$ where $(y_1,\dots,y_h\}=\n.$ Then we have $$I=(\underset{\underset{(i,j)\not=(1,2)}{1\le i < j \le h}}{y_iy_j},
\underset{3\le j \le h}{y_j^2-y_1^s},
y_2^2-y_1^{r+1}y_2-w( \underline{y})y_1^{s-t+1},y_1^ty_2).$$
  By Lemma \ref{y2} we have $y_2=cx_1^{s-t}+d$ with $c\in R$ and $d\in (x_2,\dots,x_h),$ so that $$y_2^2=c^2x_1^{2(s-t)}+2cdx_1^{s-t}+d^2.$$ Since $s\ge 2t,$ we have $s-t\ge t$ which implies $dx_1^{s-t}\in I.$ Since $t\ge r+2,$  we get also $$2(s-t)\ge 2t\ge 2(r+2)\ge r+3$$ which implies $x_1^{2(s-t)}\in \n^{r+3}.$ Finally, since $d\in (x_2,\dots,x_h)$ and $x_j^2\in I+\n^s$ for every $j\ge 3,$ we get $$d^2\in (x_2^2,\dots,x_h^2,\underset{2\le i<j \le h}{x_ix_j})\subseteq (x_2^2)+I+\n^s.$$ Now we remark that $x_2^2\in I+\n^{p+2}+\n^{s-t+1};$ since $s\ge 2t,$  $t+1\ge r+3$  and $p\ge r+1,$ this implies $x_2^2\in I+\n^{r+3}.$ Thus we get   $d^2\in I+\n^{r+3}+\n^s\subseteq I+\n^{r+3}$ because $r+3\le t+1\le s.$

  Putting all toghether, we get $y_2^2\in I+\n^{r+3};$ from this we get $$y_1^{r+1}y_2+w( \underline{y})y_1^{s-t+1}\in I+\n^{r+3}$$ which implies
  $y_1^{r+1}y_2\in I+\n^{r+3}$ because $s-t+1\ge r+3. $ This is a contradiction because $r+2\le t.$ \end{proof}

\section{The case (s,t) is not regular} In this section we are dealing with the case when there exist an integer $r$ such that $0\le r\le t-2$ and $2(r+1)=s-t+1.$ If this is the case, we say that the couple $(s,t)$ is not regular. We have seen in Section 2 that when $(s,t)$ is not regular  we have the sporadic models $I_{0,1},I_{1,1},\dots, I_{r-1,1},I_{r+1,1}\dots,I_{t-1,1}$   and the  one dimensional families  $\{I_{r,c}\}_{c\in k^*}, \{I_{r,c+x_1}\}_{c\in k^*},\dots,\{I_{r,c+x_1^{t-r-2}}\}_{c\in k^*}.$

We are going to prove that if $s\ge 2t$ all the above ideals are pairwise non isomorphic. We first remark that, if $s\ge 2t$ and $2(r+1)=s-t+1,$ then $2(r+1)\ge t+1$ so that $2r\ge t-1\ge 1.$
Hence in this section the integer $r$ is a positive integer such that $1\le r \le t-2$ and $2(r+1)=s-t+1.$

We need the following easy remarks concerning properties of the ideal $I_{r,w}$ when $2(r+1)=s-t+1$ and  $s\ge 2t.$

\begin{lemma} Let $s\ge 2t,$  $1\le r\le t-2,$  $2(r+1)=s-t+1$ and $w$ a unit in $R.$ Further let $\cong$ denote congruence modulo the ideal $I_{r,w}+\n^{3r+2}.$

a) Every  monomial of degree $r+2$ in the variables $x_1,\dots,x_h$ is in $I_{r,w}+\n^{3r+2},$ except possibly for $x_1^{r+2}, x_1^{r+1}x_2$.

\vskip 2mm

b) $\left (\sum_{i=1}^ha_ix_i\right )^{r+1}\left (\sum_{i=2}^hb_ix_i\right )\cong a_1^{r+1}x_1^{r+1}b_2x_2.$

\vskip 2mm

c) $\left (\sum_{i=2}^hb_ix_i\right )^2\cong b_2^2x_2^2.$

\vskip 2mm

d) $\left (\sum_{i=1}^ha_ix_i\right )^{2(r+1)}\cong a_1^{2(r+1)}x_1^{2(r+1)}.$

\vskip 2mm

e) If $ax_1^{r+1}x_2+bx_1^{s-t+1}\cong 0,$ then $a\in (x_1^{t-r-1})+(x_2,\dots,x_h)$ and $b\in (x_1^r)+(x_2,\dots,x_h).$
\end{lemma}
\begin{proof} For simplicity  we denote by $I$  the ideal $I_{r,w}.$ First we prove a).  If $j\ge 3$   we have $x_j^2\in I+\n^s\subseteq I+\n^{3r+2}$ because,  since $t\ge r+1,$  we have $s=2(r+1)+t-1\ge 3r+2.$ On the other hand, $x_ix_j\in I$ for $1\le i<j\le h$ and $(i,j)\not= (1,2).$ Hence we need only to consider  the monomials $x_1^rx_2^2,x_1^{r-1}x_2^3,\dots, x_2^{r+2}.$ For every $j=0,\dots,r$   we have $$x_1^{r-j}x_2^{j+2}\cong x_1^{r-j}x_2^{j}(x_1^{r+1}x_2+wx_1^{s-t+1})=x_1^{2r+1-j}x_2^{j+1}+wx_1^{s-t+1+r-j}x_2^j\cong 0$$ because the second addendum is a monomial of degree $s-t+1+r=3r+2,$ the first  is a monomial of degree $2(r+1)=s-t+1\ge 2t -t+1=t+1$ and as such, by Lemma \ref{monomi},  is in $\n^s\subseteq \n^{3r+2}.$

It is clear that b) is an easy  consequence of a), while c) is trivial. So,  let us prove d). We have $2(r+1)=s-t+1\ge t+1,$ hence, by Lemma \ref{monomi}, all the addenda in $\left (\sum_{i=1}^ha_ix_i\right )^{2(r+1)}$  , except possibly for $x_1^{2(r+1)}$,  are in $I+\n^s \subseteq I+\n^{3r+2}.$

We finally  prove e).  Let's write   $a=cx_1^j+d$ with $d\in (x_2,\dots,x_h)$ and $ 0\le j\le t-r-2.$ Then we have $$(cx_1^j+d)x_1^{r+1}x_2+bx_1^{s-t+1}=cx_1^{r+1+j}x_2+dx_1^{r+1}x_2+bx_1^{s-t+1}\cong 0.$$ Now,  by part a),  we have \begin{equation}\label{kk}dx_1^{r+1}x_2=(d_2x_2+\dots+d_hx_h)x_1^{r+1}x_2\cong d_2x_1^{r+1}x_2^2\cong 0.\end{equation} From this it follows that $$cx_1^{r+1+j}x_2\in I+\n^{3r+2}+\n^{s-t+1}\subseteq I+\n^{r+j+3} $$ because $$r+j+3\le r+t-r-2+3=t+1\le s-t+1=2r+2\le 3r+2.$$ Since $r+j+2\le t,$ by (\ref{base}) this implies $c\in \n,$ so that we can write $a=cx_1^{j+1}+d$ with $d\in (x_2,\dots,x_h).$ Going on in this way, clearly we get the conclusion $a\in (x_1^{t-r-1})+(x_2,\dots,x_h).$

This enables us to write $a=fx_1^{t-r-1}+g$ with $g\in (x_2,\dots,x_h).$  Then we have $$ax_1^{r+1}x_2=(fx_1^{t-r-1}+g)x_1^{r+1}x_2=fx_1^tx_2+x_1^{r+1}x_2g\cong 0,$$ because, as in (\ref{kk}), $x_1^{r+1}x_2g\in I+\n^{3r+2}.$
The assumption becames now $bx_1^{s-t+1}\cong 0.$

Let's write  $b=hx_1^j+k$ with $0\le j\le r-1$ and $k\in (x_2,\dots,x_h).$ Since  $s-t+1\ge t,$ we have $kx_1^{s-t+1}\in I,$  which implies $$0\cong bx_1^{s-t+1}\cong  (hx_1^j+k)x_1^{s-t+1}\cong hx_1^{s-t+1+j}.$$ This means $hx_1^{s-t+1+j}\in I+\n^{3r+2};$ now, since $r\ge j+1,$ we have
 $$3r+2\ge 2r+j+3=(s-t+1+j)+1.$$ Hence $hx_1^{s-t+1+j}\in I+\n^{3r+2}\subseteq I+\n^{(s-t+1+j)+1}.$ Since $s-t+1+j\le s-t+1+r-1=s-t+r\le s,$ by (\ref{base}) this implies $h\in \n.$ As before, going on in this way, we get the conclusion
$b\in (x_1^r)+(x_2,\dots,x_h).$
\end{proof}

We are ready to prove the main result of this section.
\begin{theorem}\label{crucial} Let $(s,t)$ be a non regular couple and $r$ the integer such that $1\le r \le t-2$ and $2(r+1)=s-t+1.$ If $s\ge 2t$ and  $I_{r,z}=\varphi_{\underline{y}}(I_{r,w})$ with $(y_1,\dots,y_h)=\n,$ then $z-w(\underline{y})\in (x_1^{t-r-1})+(x_2,\dots,x_h).$
\end{theorem}
\begin{proof} We have $$I:=I_{r,z}=\varphi_{\underline{y}}(I_{r,w})=(\underset{\underset{(i,j)\not=(1,2)}{1\le i < j \le h}}{y_iy_j},
\underset{3\le j \le h}{y_j^2-y_1^s},
y_2^2-y_1^{r+1}y_2-w( \underline{y})y_1^{s-t+1},y_1^ty_2).$$ By Lemma \ref{y2} we have $y_2=ex_1^{s-t}+f$ with $f\in (x_2,\dots,x_h).$ On the other hand,      for every $j=1,\dots,h$ we have $y_j=\sum_{i=1}^ha_{ji}x_i;$ since $\n=(y_1,\dots,y_h),$ the determinant of the matrix $(a_{ij})$ must be a unit in $R.$

Now we have $y_2^2\in I+\n^{r+2}+\n^{s-t+1}\subseteq I+\n^3$ because, since $r\ge 1,$ we have $r+2\ge 3$ and $s-t+1=2r+2\ge 4.$ Also, for every $j\ge 3,$ we have $y_j^2\in I+\n^s\subseteq I+\n^3,$ because $s\ge 3.$ By Lemma \ref{l2} we have    $a_{j1}\in \n$ for every $j\ge 2,$ so that $a_{11}\notin \n.$
Hence, we can simply write:
$$y_1=ax_1+f,\ \ \ a\notin \n, \ \ \ f\in (x_2,\dots,x_h).$$  We have
$$
y_2^2-y_1^{r+1}y_2-w( \underline{y})y_1^{s-t+1}=(ex_1^{s-t}+f)^2-y_1^{r+1}(ex_1^{s-t}+f)-w( \underline{y})y_1^{s-t+1}=$$ $$=e^2x_1^{2(s-t)}+2efx_1^{s-t}+f^2-ey_1^{r+1}x_1^{s-t}-y_1^{r+1}f-w( \underline{y})y_1^{s-t+1}\in I.$$
Now $$2(s-t)=2(s-t+1)-2=4(r+1)-2=4r+2\ge 3r+2,$$ so that $x_1^{2(s-t)}\in \n^{3r+2}.$ Also, since $s-t\ge t$ and $f\in (x_2,\dots,x_h),$ we have  $fx_1^{s-t}\in I.$ Finally $s-t+r+1=2r+2+r=3r+2,$ so that $y_1^{r+1}x_1^{s-t}\in \n^{3r+2}.$
This implies  that $$f^2-y_1^{r+1}f-w( \underline{y})y_1^{s-t+1}\in I+\n^{3r+2}.$$ Let's write $f=\sum_{j\ge 2}b_jx_j;$
using b) in the above Lemma we get $$y_1^{r+1}f\cong a^{r+1}x_1^{r+1}b_2x_2 \ \ \mod I+\n^{3r+2}.$$ Using c) and d) we further get $$y_1^{s-t+1}=y_1^{2(r+1)}\cong a^{2(r+1)}x_1^{2(r+1)} \ \ \mod I+\n^{3r+2},$$
$$f^2\cong b_2^2x_2^2 \ \ \mod I+\n^{3r+2}.$$ It follows that $$b_2^2x_2^2-a^{r+1}x_1^{r+1}b_2x_2-w( \underline{y})
a^{2(r+1)}x_1^{2(r+1)}\in I+\n^{3r+2}.$$ From this we get $$b_2^2(x_1^{r+1}x_2+zx_1^{s-t+1})-a^{r+1}x_1^{r+1}b_2x_2 -w( \underline{y})
a^{2(r+1)}x_1^{2(r+1)}\in I+\n^{3r+2}$$ and finally $$x_1^{r+1}x_2(b_2^2-b_2a^{r+1})+x_1^{s-t+1}(zb_2^2-w( \underline{y})a^{2(r+1)})\in I+\n^{3r+2}.$$ By e) in the above Lemma, this implies $$\begin{cases}
   b_2^2-b_2a^{r+1}\in (x_1^{t-r-1})+(x_2,\dots,x_h)   & \\
 zb_2^2-w( \underline{y})a^{2(r+1)}\in (x_1^{r})+(x_2,\dots,x_h).&
\end{cases}$$ Now, if $b_2$ would be  in $\n$, then $w( \underline{y})a^{2(r+1)}\in \n$, a contradiction because $a, w( \underline{y})\notin \n.$ Hence $b_2\notin \n$, so that  $$\begin{cases} b_2-a^{r+1}\in (x_1^{t-r-1})+(x_2,\dots,x_h)   &\\ z(b_2-a^{r+1})(b_2+a^{r+1})+a^{2(r+1)}(z-w( \underline{y}))\in (x_1^{r})+(x_2,\dots,x_h),&\end{cases}$$which implies $$z-w( \underline{y})\in (x_1^{t-r-1})+(x_1^{r})+(x_2,\dots,x_h).$$Since $2r+1=s-t\ge t,$ we have $r\ge t-r-1$ and we get $$z-w( \underline{y})\in  (x_1^{t-r-1})+(x_2,\dots,x_h),$$ as wanted.
\end{proof}

 We are ready now to state and prove the classification result for Gorenstein almost stretched Artinian algebras of type $(s,t)$ with the assumption $s\ge 2t.$
 \begin{theorem} \label{finale} Let $I$ be an ideal in $R=k[[x_1,\dots,x_h]]$ such that $R/I$ is almost stretched and Gorenstein of type $(s,t)$ with $s\ge 2t.$

If $(s,t)$ is regular, then $I$ is isomorphic to one and only one of the following ideals:
$$I_{0,1},\  I_{1,1},\dots , I_{t-1,1}$$

If  $(s,t)$ is not regular and $r$ is the integer such that $2(r+1)=s-t+1,$ then $I$ is isomorphic to one and only one of the following ideals:

  $$I_{0,1},\dots, I_{r-1,1},\{I_{r,c}\}_{c\in k^*},\{I_{r,c+x_1}\}_{c\in k^*},\dots,\{I_{r,c+x_1^{t-r-2}}\}_{c\in k^*},I_{r+1,1},\dots , I_{t-1,1}$$
\end{theorem}
\begin{proof} By Theorem \ref{modelli} we need only to prove the "only one" part of the statements.  By Proposition \ref{I0},  the ideal $I_{0,1}$ is not isomorphic to any of the other ideals. By Theorem \ref{noniso}, the ideals  $I_{1,1},\  I_{2,1},\dots , I_{t-1,1}$ are pairwise non isomorphic. This proves the result in the case the couple $(s,t)$ is regular.

Let us assume now that $(s,t)$ is not regular and  let $r$ be the integer such that $1\le r\le t-2$ and $2(r+1)=s-t+1.$
By Theorem \ref{noniso},  we need only to prove that however we choose two ideals in the list $\{I_{r,c}\}_{c\in k^*},\{I_{r,c+x_1}\}_{c\in k^*},\dots,\{I_{r,c+x_1^{t-r-2}}\}_{c\in k^*},$ they  are never isomorphic.

 Let
$z:=c+x_1^i,$ $w:=d+x_1^j$ with $c,d\in k^*$ and $1\le i\le j\le t-r-2.$ By contradiction, let us assume that $I_{r,z}\simeq I_{r,w};$ this means that $I_{r,z}=\varphi_{\overline{y}}(I_{r,w})$ with $(y_1,\dots,y_h)=\n.$ By Theorem \ref{crucial} we get $$z-w(\overline{y})=c+x_1^i-(d+y_1^j)\in (x_1^{t-r-1})+(x_2,\dots,x_h).$$ Since $t\ge r+2,$ we get $c-d\in \n$ which implies  $c=d,$ because  $c,d\in k^*.$  Thus we get $x_1^i-y_1^j\in (x_1^{t-r-1})+(x_2,\dots,x_h)$
which implies $i=j$ because $i\le j\le t-r-2.$ Hence $z=w$ as required.

In the same way we can prove that $I_{r,c}\simeq I_{r,d}$ implies $c=d$ and $I_{r,c}\simeq I_{r,d+x_1^j}$ implies
$j=0.$ The conclusion follows also in the case $(s,t)$ is not a regular couple.

\end{proof}

\section{Examples and remarks}

\vskip 2mm 1. Let us look at the Hilbert function $\{1,3,2,2,2,1,1,1,1\}$; this is  of type $s=8, t=4$ with $h=3.$ Since $s-t+1=5$ the couple $(8,4)$ is regular. The isomorphism classes  are represented by the following ideals $$(x_1x_3,x_2x_3,x_3^2-x_1^8, x_2^2-x_1^{n+1}x_2-x_1^5,x_1^4y)$$ for $n=0,1,2,3.$  Hence we have a finite number of isomorphism classes.

\vskip 2mm \noindent 2. If we consider the Hilbert function $\{1,2,2,2,1,1,1,\}$ then $h=2,$ $t=3$ and $s=6=2t.$ We have $s-t+1=4=2(1+1)$
so that  the couple $(6,3)$ is not regular.  The isomorphism classes are represented by the following ideals
$$I_{0,1}=(y^2-xy-x^4,x^3y),\ \  I_{2,1}=(y^2-x^3y-x^4,x^3y)$$ and $$\ \{I_{1,c}\}_{c\in k^*}=\{(y^2-x^2y-cx^4,x^3y)\}_{c\in k^*}.$$ This example has been studied in \cite{EV} with different methods. It is the first case where an infinite number of isomorphism classes  arises, namely two sporadic models plus a one dimensional family. The understanding of this difficult example was the motivation of this work. Notice that the length is 

\vskip 2mm \noindent 3. We can produce examples where there are several one dimensional families of models. Take the Hilbert function $\{1,2,2,2,2,2,1,1,1,1,1\};$ then $h=2,$ $t=5,$  $s=10$ and  the couple $(10,5)$ is not regular. The isomorphism classes are represented by the following ideals:
$$\begin{cases}
    (y^2-x^{r+1}y-x^6,x^5y)  & \text{for} \ \ \ r=0,1,3,4\\
     (y^2-x^3y-cx^6,x^5y) & \text{for}  \ \ \ c\in k^*\\
   (y^2-x^3y-cx^6-x^7,x^5y)  & \text{for} \ \ \ c\in k^*.\end{cases}$$

   \vskip 2mm \noindent 4. The above description of the isomorphism classes of almost stretched Gorenstein algebras with a given Hilbert function is no more available if we do not assume $s\ge 2t.$ For example let $t=3,$ $s=5$ and $h=2,$ corresponding to the Hilbert function $\{1,2,2,2,1,1\}.$ Then $s=2t-1$ and the  couple $(5,3)$ is regular; nevertheless we can prove that $I_{2,1}\sim I_{1,1},$ thus contradicting the conclusion of Theorem \ref{noniso}.

   We have $I:=I_{1,1}=(y^2-x^2y-x^3,x^3y)$ and $I_{2,1}=(y^2-x^3,x^3y).$ Let us change the variable as follows:
   $$\begin{cases}
  z=9x+y    & \text{ } \\
  w=-27y+xy+9x^2.   & \text{}
\end{cases}$$ We have $\n=(z,w)$ and   the following congruences mod $I$ hold true:
$$xy^2\cong x^4, \ \ x^2y^2\cong x^5, \ \  y^3\cong x^5,  \ \ xy^3\cong 0, \ \ y^4\cong (x^2y+x^3)^2\cong 0.$$ From this we get $$w^3=(-27y+xy+9x^2)^3\simeq (-27)^3x^5+3(-27)^29x^5\simeq 0$$ and $$w^2-z^3=(-27y+xy+9x^2)^2-(9x+y)^3\simeq 0.$$ This proves that $I\supseteq (w^3,w^2-z^3)=(w^2-z^3,wz^3);$ by computing the Hilbert function of the last ideal we can see that $$I=I_{1,1}= (w^3,w^2-z^3)=\phi_{\{z,w\}}(I_{2,1}),$$ as claimed.

As a consequence we get that the family of ideals $I$ such that $R/I$ is Gorenstein with Hilbert fuction $\{1,2,2,2,1,1\}$ has two isomorphic classes, those corresponding to the following ideals: $$ (y^2-xy-x^3,x^3y)\sim(xy,y^4-x^5),$$ $$ (y^3,y^2-x^3)=(x^3y,y^2-x^3)\sim (y^2-x^2y-x^3,x^3y).$$ We need to remark that the isomorphism $(y^2-xy-x^3,x^3y)\sim(xy,y^4-x^5)$ comes from the following easy claim: $$(x-y+xy)(y+x^2+x^3)\in (y^2-xy-x^3,x^3y).$$

\vskip 2mm \noindent 5. The last remark is dealing with the case $s=t+1,$ which is far away from the basic  assumption $s\ge 2t$ which we  used in the paper. The couple $(t+1,t)$ is not regular,  the critical value beeing $r=0.$ We are able to prove that every  ideals $I$ such that $R/I$ is Gorenstein with Hilbert fuction $\{1,h,2,2,\dots,2,1\}$ is isomorphic to one of the following ideals:$$I_{t-1,1}=(\underset{\underset{(i,j)\not=(1,2)}{1\le i < j \le h}}{x_ix_j},
\underset{3\le j \le h}{x_j^2-x_1^s},
x_2^2-x_1^2,x_1^tx_2)$$ $$I_n:=(\underset{\underset{(i,j)\not=(1,2)}{1\le i < j \le h}}{x_ix_j},
\underset{3\le j \le h}{x_j^2-x_1^s},
x_2^2-x_1^n,x_1^{t+1}-x_1^tx_2),  \ \ \ n=3,\dots,t$$
$$(\underset{\underset{(i,j)\not=(1,2)}{1\le i < j \le h}}{x_ix_j},
\underset{3\le j \le h}{x_j^2-x_1^s},
x_2^2,x_1^{t+1}-x_1^tx_2).$$ If $t\ge 4,$ we  are not able to prove that these ideals are pairwise non isomorphic.

\vskip 3mm The last two examples show that the problem of the classification up to ismorphism of almost stretched Gorenstein algebras with a given Hilbert function becomes more difficult  when $s\le 2t-1.$ This is the reason why, at the moment, we really need the assumption  $s\ge 2t.$

\providecommand{\bysame}{\leavevmode\hbox to3em{\hrulefill}\thinspace}

\end{document}
\end